\documentstyle[12pt]{article}
\begin{document}
\def\A{{\cal A}}
\def\B{{\cal B}}
\def\C{{\cal C}}
\def\D{{\cal D}}
\def\E{{\cal E}}
\def\F{{\cal F}}
\def\G{{\cal G}}
\def\I{{\cal I}}
\def\J{{\cal J}}
\def\K{{\cal K}}
\def\M{{\cal M}}
\def\N{{\cal N}}
\def\Q{{\cal Q}}
\def\R{{\cal R}}
\def\T{{\cal T}}
\def\U{{\cal U}}
\def\V{{\cal V}}
\def\W{{\cal W}}
\def\X{{\cal X}}
\def\Y{{\cal Y}}
\def\Z{{\cal Z}}
\def\be{\begin{equation}}
\def\ee{\end{equation}}
\def\la{\label}
\def\sc{\scriptstyle}
\def\tf{\tilde\F}
\def\da{\dagger}
\def\bC{{\bf C}}
\def\bG{{\bf G}}
\def\bI{{\bf I}}
\def\bM{{\bf M}}
\def\bP{{\bf P}}

\begin{center}{\Large Application of Tree-like
Structure of Graph to Matrix Analysis.}
\end{center}

\centerline{\bf Buslov V.A.}
\vskip4mm

\begin{center}{\it Abstract}\end{center}

{\small Formulas for matrix determinants, algebraic adjunctions,
characteristic poly\-no\-mial coefficients, components of eigenvectors are
obtained in the form of signless sums of matrix elements products taking by
special graphs.  Signless formulas are very important for singular and
stochastic problems.  They are also useful for spectral analysis of large
very sparse matrices.}

\section{Introduction}

Graph theory is a natural instrument for matrix determinants and combined
features calculating.  It is clear that determinants can be expressed in
form of alternating sum of matrix elements products taking by some graphs
corresponding to index permutation, each term sign is determined by
permutation even. However, the alternation itself is an essential
obstruction both in numerical calculating and in analytic research of
spectrum properties. It concerns for example cases when matrix elements
depend on small (large) parameter and possess on it different orders.
Thus under singular perturbed equations research \cite{L} (in particular
such as Fokker-Plank equations) necessity appears for spectral analysis of
big order matrix having exponentially small elements \cite{V1,V4}.  The
same situation appear under consideration of Markov's chains connected
with diffusion process \cite{VF} where the questions of stochastic
continuity, possible subprocess kinds and their structure properties run
into the necessity of spectral analysis just such kind of matrixes.

Our aim is to obtain the formulas for matrix determinants, algebraic
adjunctions, characteristic polynomial coefficients, compo\-nents of
eigen-vectors in the form of signless sums of matrix elements products
taking by special graphs. It turns out that one can get such formulas in
terms of so called "tree" structure of some graph corresponding to matrix.

The first step in this direction was made by Kirchhoff \cite{K} who
computed the number of connected subgraphs containing all vertices and
containing no circuits (spanning trees). This number turned out to be
equal to cofactor of any element of so called conductivity or Kirchhoff
matrix of non-directed graph.  Later this theorem was generalized for all
coefficients of characteristic polynomial of this matrix, and also for
cases of directed graphs where every arc (ordered pairs of vertices)
possesses some quantity called weight \cite{FS}. Here we get the same type
formulas not for special matrices (as Kirchhoff type) but for arbitrary
ones.

\section{Main definitions and designations}

Unification of designations and even terminology proper is not complete
yet in graph theory. So firstly we adduce the necessary definitions and
notations.

Let $\ G\ $ be digraph (directed graph). We use  $\ \V G\ $ and
$\ \A G\ $  to denote the set of vertices and arcs of $\ G\ $. The
subgraph $H$ of $G$ is called factor if $\V H=\V G$. The outdegree
(indegree) of the vertex $\ i\ $ (the number of arcs going out of (into)
$\ i\ $) we denote $\ d^+(i)\ $ ($d^-(i)$).

If $G$ is digraph in which every arc has its own weight  $g_{ij}$
(weighted adjacencies digraph), corresponding matrix $\bG = \{g_{ij}\}_
{i,j=1}^{N}$ is called generalized adjacency matrix (the element
$g_{ij}=0$ if there is no arc $(i,j)$ in $G$).

The sequence of following each other arcs along their orientations is
called a way if all vertices besides possibly the uttermost ones are
different. The way connecting the vertices $m$ and $n$ we denote
$m\cdot n$. The cyclic way is called dicircuit. Linear digraph is
digraph every vertex of which has unit in- and out-degree. So it
consists of dicircuits.

We associate with any weighed adjacencies digraph $G$ the quantity
$\pi_G$ by the rule

$$\pi_G=\prod_{(i,j)\in\A G}g_{ij} \ $$
which is naturally to call by productivity.

Later on forests are the main graph theory object we use. As known there
are two forest kinds in digraph situation. Here we call by forest
digraph without dicircuits in which every vertex outdegree is equal to 0
or 1 ($d^+(i)=0,1$). The only vertex of tree (component of forest)
having zero outdegree ($d^+(i)=0$) we call root.

Let $\N $ be finite set and $\W $ is some subset of $\N $.
By $\F ^k_\W $, $k\le|\N |-|\W |$ we denote the set of
forests $F$ obeying the following conditions

1) ${\cal V}F={\cal N}$;

2) $F$ consists of exactly $k+|\W |$ trees;

3) The set of roots of $F$ contain $\W $ as a subset.
\\
Suppose also $\F ^k = \F ^k_{\emptyset}$, $\F _\W  = \F ^0_\W $. Note that
the set  $\F ^0$ is empty set and $\F ^{|\N  |}$ consists of the only
empty forest having only roots and no arcs.

If some additional condition $z$ is put on forests from $\F ^k_
\U $ we denote such set of forests as $\F ^k_\U (z)$. Other
necessary utilized notations we sign as necessary and also we use
sometimes the  term "graph" in wide sense designating by it digraphs with
weighted adjacencies too.

\section{Circuit- and tree-like structures}

The matrix (graph) spectral analysis can be carried out using some form of
characteristic polynomial of matrix itself or some special matrixes
(graphs) constructed from it. Thus it is valid known "theorem on
coefficients for digraphs" \cite{CDZ}.

{\bf Theorem}. {\em Let

$$
\det (\lambda\bI -{\bf A})=\lambda^N+a_1\lambda^{N-1}+\cdots +a_N
$$
be characteristic polynomial of arbitrary digraph $A$ with weighted
adjacencies $a_{ij}$. Then

\be
a_i=\sum_{L\in {\cal L}_i}(-1)^{p(L)}\pi_L \ , \ \
i=1,2,\cdots ,N \ ,
\la{circuit}
\ee
where ${\cal L}_i$ is the set of all linear directed subgraphs $L$ of
graph $A$ with exactly $i$ vertices; $p(L)$ means the number of
components (dicircuits) of $L$}.

Coefficients in (\ref{circuit}) are expressed in "circuit" structure
of $A$, and this theorem is not more than rephrasing from the standard
determinant notation $|\lambda\bI -{\bf A}|$ in the form of matrix
elements products sum with sign determined by substitution even (the
number of dicircuits) into graph terms.

In terms of "tree" structure it is known the characteristic polynomial
expression not for matrix $\bG $ itself but for its Kirchhoff matrix
$\bC $ (or conductivity matrix) determined like

$$ \bC ={\bf D}-\bG  \ ,  $$
where ${\bf D}$ is weighted powers matrix

$$ {\bf D}\equiv {\rm diag}(\sum_{j=1}^Ng_{1j},\sum_{j=1}^Ng_{2j},\cdots ,
\sum_{j=1}^Ng_{Nj}  ) \ .$$
Corresponding expression has a form \cite{FS}:

\be
\det (\lambda\bI -\bC )=(-1)^N\sum\limits_{k=0}^N(-\lambda
)^k \Bigl[\sum\limits_{F\in\F ^k(G)}\pi_F\Bigr] \ ,
\label{ein}
\ee

$$ \pi_F=\prod\limits_{(i,j)\in\A F}f_{ij}=
\prod\limits_{(i,j)\in\A F}g_{ij} \ .
$$
Here the set of forests containing directly $k$ trees and being factors
of $G$ is designated by $\F ^k(G)$.  Note that as the sum of elements
along every line of $\bC $ is equal to zero, so its determinant is equal
to zero too and the sum in (ref{ein}) one can lead from $k=1$.  In the
following in clear cases we  omit indication on graph.

\section{The characteristic polynomial in tree-like structure terms}

Knowing the characteristic polynomial expression of the admittance matrix
$\bC $ in terms of tree structure it is not hard to get analogous
expression for characteristic polynomial of the matrix $\bG $ itself.
For this aim let us construct from $G$ some new graph $G^\da$ by
the next rule. Let $\N ^\da =\N \cup\{ \da \}$ be the
set of vertices $\N =\{ 1,2,\cdots , N\}$ of $G$ which is
supplemented by some new vertex  designated $\da $. Let also add to
$G$ the vertex $\da $ and lay on out of every vertex $i\in\N $
the arc $(i,\da )$ with weight $g_{i\da}=-\sum\limits_{j=1}^N
g_{ij}$, and remove all loops $(i,i)$.  We denote the obtained graph by
$G^\da $ ( $G^\da :  \V G^\da =\N ^\da$,
$\A G^\da = \{ \A G\setminus\mathop{\cup}\limits_{i\in
\N }(i,i)\} \mathop{\cup}\limits_{i\in \N }(i,\da )$ \ , weights of
arcs are equal $g_{ij} $ , $i\in \N $ , $j\in \N ^ \da $).  Corresponding
generalized adjacencies matrix $\bG ^\da $ of $G^\da $ has the following
form

$$ \mbox{\bf\Large G$^\da $}=\pmatrix{
                 0 & g_{12} & g_{13} & \ldots &g_{1N}& g_{1\da }\cr
                 g_{21} & 0 & g_{23} & \ldots &g_{2N}& g_{2\da } \cr
                 g_{31} & g_{32} & 0 & \ldots &g_{3N}& g_{3\da } \cr
            \vdots & \vdots & \vdots &\ddots  & \vdots &\vdots \cr
                 g_{N1} & g_{N2} & g_{N3} & \ldots & 0 &g_{N\da } \cr
		 0     &  0     &    0   &  \ldots & 0 & 0 \cr
	      } \ .
$$
Note that, if $\bG =\bP -\bI $, where $\bP $ is  probability
matrix setting finite Markov's chain with killing, so the
quantity $g_{i\da }$ is the probability of killing of the process
if it is in state $i$. This has a sense of probability of outcoming to
the bounder, so the additional vertex $\da $ could be
interpreted in some sense as a bounder of the finite set $\N =\{
1,2,\cdots ,N\}$.

It is easy to see, that $(N+1)\times (N+1)$ admittance matrix ${\bf
C}^\da $ of graph $G^\da $ has the form

$$ \mbox{\bf\Large C$^\da $} = - \left(
\array{c|c} \mbox{\bf\Large G} & \array{c} g_{1\da }\\ g_{2\da }  \\
\vdots \\ g_{N\da } \endarray \\ \hline \array{cccc} 0 & 0 & \cdots & 0
\endarray
&
0
\endarray
\right) \ ,
$$
so, using (\ref{ein}) we get the chain of equations

$$
\det (\lambda\bI -\bG )={1\over\lambda}
\det (\lambda\bI +{\bf C^\da })
={1\over\lambda}\sum\limits_{k=1}^{N+1}
\lambda^k\Bigl[\sum\limits_{F\in\F ^k(G^\da )}\pi_F\Bigr] =
$$

$$
=\sum\limits_{k=0}^{N}
\lambda^k\Bigl[\sum\limits_{F\in\F ^{k+1}(G^\da )}\pi_F\Bigr]
\ .
$$
Note, that since the vertex $\da $ in graph $G^\da $ has zero
outdegree and hence it is a root in every forest $F\in\F^{l}(G^\da
)$, so the sets  $\F^{k+1}(G^\da )$ and $\F ^{k}_\da (G^\da )$ coincide.
Thus it is valid \\

{\bf Theorem 1.} {\em Characteristic polynomial of an
arbitrary $N\times N$ matrix $\bG$ can be expressed in the form

\be
\det
(\lambda\bI -\bG )=\sum\limits_{k=0}^{N}
\lambda^k\Bigl[\sum\limits_{F\in\F ^{k}_\da }\pi_F\Bigr]
\ ,\ \ \pi_F=\prod\limits_{(i,j)\in\A F}g_{ij}\ ,
\label{char}
\ee
where $\F ^k_\da \equiv \F ^k_\da (G^\da ) \ .$}

Under $\lambda =0$ we obtain obvious \\
{\bf Consequence.} {\em The determinant of  $N\times N$ matrix $\bG$  can
be expressed in the form

\be
\det\bG =(-1)^N\sum\limits_{F\in\F _\da }\pi_F  \ ,
\label{det}
\ee

$$ \F _\da \equiv\F _\da (G^\da ) \ .$$}

Let  $\R $ be a subset of $\N $. Designate by $\bG _{\cal RR}$ matrix
obtained from  $\bG $ by striking out the columns and lines with numbers
$i\in \R $.  So, $\bG _{\cal RR}$ is a diagonal minor of $\bG $ of $(N-|\R
|)$-th order. The corresponding to it digraph we denote by $G_{\cal RR}$.

{\bf Consequence of consequence.} {\em The determinant of minor $\bG
_{\cal RR}$ can be expressed in the form }

\be
\det\bG _{\cal RR}=(-1)^{N-|\R |} \sum\limits_{F\in\F _{\{
\da \}\cup \R }}\pi_F \ ,
\label{minor}
\ee

$${\cal
F}_{\{ \da \}\cup\R }=\F _{\{ \da \}\cup\R }
(G^\da )\ .
$$

{\bf Proof.} By formula (\ref{det})

$$ \det{\bf
G}_{\cal RR}=(-1)^{N-|\R |}\sum\limits_{F\in\F _\da
(G^\da _\R )} \pi_F \ .  $$
 Let us keep in $\bG _{\cal RR}$ the same numeration of elements as it is
in matrix $\bG $. So the elements $(g_{\cal RR})_{il}$ of the
corresponding matrix $\bG _{\cal RR}^ \da $ (weights of arcs $(i,l)$ of
graph $G_{\cal RR}^\da $) are equal:

 $(g_{\cal RR})_{il}=g_{il} \  , \  i,l\in \N \setminus
 \R  \ , $

  $(g_{\cal RR})_{il}=0 \ \, \{ i,l\}\cap
 \R \not= \emptyset , $

 $(g_{\cal RR})_{i\da }=\sum\limits_{m\in \{ \da \}\cup\R }
g_{im} \ , i\in \N \setminus \R  \ . $ \\ 
So productivity $\pi_F$ of any forest $F\in\F _\da(G^\da _{\cal RR})
$ represents the productivity $\pi_{F'}$ of some forest of the set  \ $\F
 _{\{ \da \}\cup \R }({ G}^\da )$.  The sum of productivities $\pi_F$
 along all forests $F\in\F _\da({ G}^\da_{\cal RR})$ exhausts the set $\F
_ {\{ \da \}\cup \R }({ G}^\da )$, which proves (\ref{minor}).  \\

\section{Formulas for components of eigenvectors}

If it is known the eigenvalue $\lambda$ of $\bG $, so to calculate
the components $v_m$ of corresponding eigenvector $\vec v$ it is
necessary to decide the standard system $(\lambda\bI -\bG )\vec v=0
\ $. Let for certainty the $n$-th component of eigenvector $\vec v$ be not
equal to zero. Without loss of generality  one can accept it be equal to
one:  $v_n=1$. Then we have for the rest of components following the
Kramer rule

\be
v_m={\det\Delta_{nm}'(\lambda) \over\det\Delta_{nn}(\lambda )}\ ,
\label{vector}
\ee
where

$$\Delta_{nn}(\lambda )=\lambda\bI -\bG _{nn} \  ,
$$
$\bG _{nn}$ is algebraic adjunct of the element $g_{nn}$ of
matrix $\bG $, and  $\Delta_{nm}'(\lambda )$ is a matrix obtained from
$\Delta_{nn}(\lambda )$ by substitution of the $m$-th column of
$\lambda\bI -\bG $ by $n$-th one with negative sign.

The expression for $\Delta_{nn}(\lambda)$ it is easy to get using
already obtained formula (\ref{minor}) concerning diagonal minors  ${\bf
G}_{{\cal RR}}$ of matrix  $\bG $. Since the coefficient at $\lambda^k$
at the expression of characteristic polynomial of algebraic adjunct
$\bG _{nn}$ is itself a sum of determinants of diagonal minors of
matrix \ $\bG $ \ of \ $(N-k-1)$-th order with sign \ $(-1)^{N-k-1}$ \
and not including $n$-th column and  $n$-th line so one can write

$$
\det (\lambda{\bf
I}-\bG _{nn})=\sum\limits_{k=0}^{N-1}\lambda^k\Bigl[\sum\limits_
{\sc \R\in\N, \  n\in\R \atop \sc
|\R |=k+1}\det ({-\bf G}_{{\cal RR}})\Bigr]= $$

$$=\sum\limits_{k=0}^{N-1}\lambda^k \Bigl[(-1)^{N-k-1}\sum\limits_
{\sc \R\in\N ,\ n\in \R \atop\sc
|\R |=k+1}\det\bG _{{\cal RR}}\Bigr]= $$

$$
=\sum\limits_{k=0}^{N-1}\lambda^k\Biggl[\sum\limits_
{\sc \R\in\N ,\ n\in\R \atop\sc
|\R |=k+1}\biggl\{\sum\limits_{F\in \F _{\{ \da \} \cup
\R }(G^\da )}\pi_F\biggl\}\Biggr] \ .$$
The sum at square brackets in the last expression obviously is equal to
the sum of productivities $\pi_F$ of all forests belonging to the
unification

$$ {\mathop \bigcup\limits_ {\sc \R\in\N ,\ n\in\R \atop\sc |\R |=k+1}} \F
_{\{ \da \} \cup\R }(G^\da ) \ , $$
which in its turn is equal to $\F_{\{ \da ,n\} }^k(G^\da )$.  So, the
denominator at (\ref{vector}) has the form

\be
\det\Delta_{nn}(\lambda )=\det (\lambda\bI -\bG _{nn})=
\sum\limits_{k=0}^{N-1}\lambda^k\Biggl[\sum\limits_
{F\in \F _{\{ \da ,n\}}^k}\pi_F\Biggr] \ . \label{znamenatel}
\ee

$$ \F _{\{ \da ,n\}}^k\equiv \F _{\{ \da ,n\} }^k(G^\da )\ , $$

It is more difficult to obtain the formula for numerator at (\ref{vector})
having the form

\be
\det\Delta_{nm}'(\lambda)=\sum\limits_{k=0}^{N-1}\lambda^k
\Biggl[\sum\limits_ {F\in \F _{\{ \da ,n\}}^k(m\cdot n)}\pi_F
\Biggr] \ ,
\label{chislitel}
\ee
where by $\F _{\{\da ,n\}}^k(m\cdot n)$ the subset of $\F _{\{ \da 
,n\}}^k( G^\da )$ having $m\cdot n$-walk is denoted. The summing over $k$ 
one can lead till $N-2$, because the set $\F _{\{ \da ,n\}}^{N-1}(m\cdot 
n)$ is empty.

It is sufficiently to establish (\ref{chislitel}) under $\lambda =0$,
because the coefficients at powers of $\lambda $ are the diagonal minors
of the same matrix and are utterly of analogical to each other form.

Let $\Delta_{nm}'\equiv\Delta_{nm}'(0)$. It is necessary to prove that

\be
\det\Delta_{nm}'=\sum\limits_ {F\in \F _{\{ \da ,n\}}(m\cdot n)}
\pi_F \ .
\label{0}
\ee
In passing let us remark that using (\ref{0}) it is easy to obtain the
expression for algebraic adjunct $G_{nm}$ of the element $g_{nm}$.
Actually, matrices $-\Delta_{nm}'$ and $G_{nm}$ differ only in such a way.
In $-\Delta_{nm}'$ (if for example $n<m$) $n$-th column of matrix $G_{nm}$
is at the $(m-1)$-th place and is of  opposite sign. Hence

$$G_{nm}=(-1)^{N+m-n-1}\sum\limits_{F\in\F _{\{\da ,n\} }}\pi_F \ .
$$

To prove (\ref{0}) let us exchange in $\bG $ places of $n$-th  and
$m$-th columns. Obtained auxiliary matrix we denote by ${\bf H}$.
Nondiagonal elements $h_{ij}$ of  ${\bf H}$ and corresponding expanded
matrix ${\bf H}^\da $ are equal

i) $h_{ij}=g_{ij} \ , \
j\not=n,m$ \ ,

ii) $ h_{in}=g_{im} \ , \ \ h_{im}=g_{in} \ $ .

It is easy to see that algebraic adjunct ${\bf H}_{nn}$ of the
element $h_{nn}$ of ${\bf H}$ differs from matrix $(-\Delta_{nm}')$ by
only the sign of $m$-th column. So taking into account (\ref{znamenatel})

\be
\det\Delta_{nm}'=-(-1)^{N-1}\det {\bf H}_{nn}=
-\sum\limits_
{F\in \F _{\{ \da ,n\}}(H)}\pi_F \ .
\label{H}
\ee
Now in (\ref{H}) it is necessary to cross from the sum of productivities
$\pi_F$ of factor forests of $H^\da $ to the sum of productivities
$\pi_F'$ of corresponding subgraphs of graph $G^\da $.  Notice
that in $H^\da $  weight of arc $(m,n)$ is equal to minus sum of
weights of arcs $(m,i)$ going out of the vertex $m$ in graph $G^\da
$:

\be
h_{mn}=g_{mm}=-\sum\limits_{i\in \N ^\da \setminus\{ m\} }g_{mi}
\ .  \label{h=-g}
\ee
That is why we divide the set of forests
$\F _{\{ \da ,n\}}( H^\da )$ in to two nonintersecting
subsets: the set $\F _{\{ \da ,n\}}( H^\da ; (m,n))$ of
forests including the arc $(m,n)$ and the set $\F _{\{ \da ,n\}}
(H^\da ;\rceil (m,n)) $ of forests without this arc.  Then

\be
\sum\limits_{F\in \F _{\{\da ,n\}}( H^\da )}\pi_F=
\sum\limits_ {F\in \F _{\{ \da ,n\}}(H^\da ;(m,n))}\pi_F
+\sum\limits_ {F\in \F _{\{ \da ,n\}}(H^\da ;\rceil
(m,n))}\pi_F \ .
\label{sum h=sum g}
\ee
The set of arcs $\A F$ of the forest $F\in \F _{\{\da
,n\}}( H^\da ;(m,n))$ is representable in a form  $\A F=(m,n) \cup\A
P$ , where $P$ is some graph belonging to the set $\F _{\{ \da ,n,m\}}(
H^\da )$. By it as $h_{ij}=g_{ij}$,  $j\not= m,n$ the sets $\F
_{\{\da ,n,m\}} (H^\da )$ and $\F _{\{ \da ,n,m\}}( G^\da )$ coincide.
Hence, taking into account (\ref{h=-g})

$$
\sum\limits_{F\in \F _{\{ \da ,n\}}(H^\da ;(m,n))}\pi_F=
-\sum\limits_{i\in \N ^\da \setminus\{ m\} }g_{mi}
\Biggl[\sum\limits_{F\in \F _{\{ \da ,n,m\}}(G^\da
)}\pi_F\Biggl] \ .
$$
Note that any forest $F\in \F _{\{\da ,n,m\}}( G^\da )$
consists of exactly of three trees  $T_\da , \ T_n \ $ and $T_m$  with
 roots correspondingly  $\da , \ n$ and $ m$.  Therefore the addition to
 such forest the arc $(m,i)$ with weight  $g_{mi}$ leads either  to graph
belonging to $\F _{\{\da ,n\}}(G^\da )$ (if $i\in\V T_
\da \cup\V T_n$), or (if $i\in\V T_m$) to graph, in which
outdegree of every vertex is equal to one with the exception of $\da $ and 
$n$:\   $d^+(j)=1 \ j\not=\da ,n$ , $d^+(\da )=d^+(n)=0$, and there is 
 exactly one circuit in graph and this circuit contains the vertex 
 $m$. We denote the set of such graphs by ${\cal O}_m(G^\da )$.  So

\be
\sum\limits_ {F\in \F _{\{ \da ,n\}}( H^\da ; (m,n))}\pi_F=
 -\sum\limits_{F\in \F _{\{ \da ,n,\}}(G^\da )}\pi_F-
 \sum\limits_ {F\in {\cal O}_m(G^\da )}\pi_F \ .
\label{middle}
\ee

Let us now consider the last sum at the right part of (\ref{sum h=sum g}).
The productivity $\pi_F$ of any forest $F\in\F _{\{\da
,n,\}}(H^\da ;\rceil (m,n))$ by force of i) ¨ ii) one can represent as
 $\ \pi_F=\pi_E$, where  $E$ is some subgraph of $G^\da $, to
obtain which from the forest $F$ one must exchange the arcs ending in $m$
into arcs ending in $n$ and back. The forest $F\in\F _{\{\da
,n,\}} (H^\da ;\rceil (m,n)) $ consists of two trees $T_\da $ and
$T_n$ with the roots $\da $ and $n$ correspondingly. By this if
$m\in\V T_\da $ so the graph $E$ is still a forest and besides
the sequence of arcs from $m$ does not lead to  $n$ (but leads to $\da
$). In the case of $m\in T_n$ in $F$  such an exchanging of arcs results
in $E$ be graph containing the only circuit and the vertex $m$ belongs
to this circuit, i.e. $E\in {\cal O}_m(G^\da )$.  Thus

\be
\sum\limits_{F\in \F _{\{ \da ,n\}}(H^\da ;\rceil (m,n))}
\pi_F= \sum\limits_ {F\in \F _{\{ \da ,n\}}(G^\da ;\rceil
 m\cdot n)}\pi_F +\sum\limits_ {F\in {\cal O}_m(G^\da )}\pi_F \ .
\label{end}
\ee
Here we denote by $\F _{\{ \da ,n\}}(G^\da ;\rceil m\cdot n)$
the subset of  $\F _{\{ \da ,n\}}(G^\da )$ not containing
$m\cdot n$-way. Combining together (\ref{H}),(\ref{sum
h=sum g}),(\ref{middle}) and (\ref{end}) we get, that products $\pi_F$
along graphs containing circuits reduce and

$$
\det\Delta_{nm}'=\sum\limits_ {F\in \F_{\{\da ,n\}}(G^\da )}\pi_F-
\sum\limits_{F\in \F_{\{ \da ,n\}} (G^\da ;\rceil m\cdot
n)}\pi_F=\sum\limits_{F\in\F_{\{ \da ,n\}}(G^\da ;m\cdot n)}\pi_F  \ ,
$$
which is the same as (\ref{0}). Thus it is proved \\
{\bf Theorem 2}. {\em Let $\lambda $ be the eigenvalue of order 1 of
matrix \bG , ${\vec v}$ | corresponding eigenvector and its $n$-th
component is not equal to zero, then accurate to constant factor the
components of $\vec v$ are following

\be
 v_n=1 \ \ , \ \ v_m={\sum\limits_{k=0}^{N-2}
\lambda^k\Bigl[\sum\limits_{F\in\F^k_{\{n,\da \}}(m\cdot n)}
\pi_F\Bigr]\over \sum\limits_{k=0}^{N-1}
\lambda^k\Bigl[\sum\limits_{F\in\F^k_{\{n,\da \}}} \pi_F
\Bigr]}\ .
\label{components}
\ee
where $\F^k_{\{n,\da \}}(m\cdot n)$ is a subset of set
 $\F^k_{ \{n,\da \}}$, consisting of forests having $m\cdot
n$-way.}

{\bf Remark.} To obtain the expression for $m$-th component of eigenvector
of the transform matrix $\bG ^{\top}$ it is necessary only to exchange 
indices $m,n$ with each other at the right part of (\ref{components}).

\section{Example}

Let us consider  an easy example demonstrating the application of singless
form technique. Let $\bM =\bP -\bI $, where $\bP $ is a
probability matrix, setting the finite Markov's chain with killing 
(that is the sum of elements along any line may be less than one, 
corresponding residual is just a killing probability), and also let the 
transition probabilities $M_{ij}=P_{ij}$, $i\ne j$, and the killing ones 
$M_{i\da }=-\sum\limits_{j=1}^NM_{ij}= 1- \sum\limits_{j=1}^NP_{ij}$ be 
exponentially small. Namely, let $N=3$ and $M_{ij}=m_{ij}e^{-V_{ij}/
\varepsilon}$ and $V_{12}$=$V_{13}=4$, $V_{1\da }=5$, $V_{21}=3$,
$V_{23}=2$, $V_{2\da }=5$, $V_{32}=1$, $V_{3\da }=4$, $V_{31}=3$;   
$\varepsilon$ is a small parameter.  Thus if  instead of diagonal elements 
$M_{ii}$ we use their expression through the quantities $M_{i\da }$ 
($M_{ii}=-M_{i\da } -\sum\limits_{j\ne i}^NM_{ij}$)  the matrix $\bM $ 
gets the form

$$\bM = \pmatrix{
-M_{12}-M_{13}-M_{1\da}& M_{12}
&M_{13} \cr M_{21}
&-M_{21}-
 M_{23}-M_{2\da} &
 M_{23} \cr M_{31}&
 M_{32} &-M_{32}
 -M_{3\da }-M_{31} \cr} \ .
$$
By virtue of lack of sign and non-negativity of transition (nondiagonal
elements of $\bM $) and killing ($M_{i\da }$) probabilities that are
used at tree-like structure formulas, it is ought to keep at the
asymptotic of characteristic polynomial coefficients (\ref{char}) only terms
reaching the maximal order on small parameter. Thus the eigenvalues
asymptotic one can extract from the equation

$$
\lambda^3+\lambda^2 m_{32}e^{-1/\varepsilon}+ \lambda
m_{32}m_{21}e^{-5/\varepsilon}+ m_{32}m_{21}m_{1\da }
e^{-10/\varepsilon}=0\ .
$$
The exponential orders $V_k=-\lim\limits_{\varepsilon\to 0} \varepsilon
\ln a_k $ of the coefficients  $a_k$ of characteristic polynomial
$\sum\lambda^k a_k$ satisfy  convex nonequalities system analogical to 
\cite{V}:  $V_{k-1}-V_k\ge V_k-V_{k+1}$, and the eigenvalues asymptotic
one can look in the form $\lambda_k^\varepsilon {\mathop{\sim}\limits_{
\varepsilon\to 0}}\Lambda_k\ e^{\frac {V_{k}-V_{k-1}}{\varepsilon}}$
\cite{V1}. Substituting model $\lambda$ in such a form we get
the following eigenvalues asymptotics:

$$\lambda_1{\mathop{\sim}\limits_{\varepsilon\to 0}}-m_{1\da }
e^{-5/\varepsilon}\  ,\ \
 \lambda_2{\mathop{\sim}\limits_{\varepsilon\to
0}}-m_{21} e^{-3/\varepsilon}\ , \ \
\lambda_3{\mathop{\sim}\limits_{\varepsilon\to 0}}-m_{32}
e^{-1/\varepsilon}\ . $$

The asymptotic of eigenvectors of the matrix $\bM $ and the transform 
matrix $\bM ^*$ one can find from (\ref{components}). Denoting by $\bC $ 
(by $\bC '$) the matrix, composed  of ultimate values of
vector-lines of $\bM $ (vector-columns of $\bM ^*$) one gets

$$\bC =\pmatrix{1&1&1\cr
                   0&1&1\cr
                   0&0&1\cr}\ ,\ \
    {\bf C'}=\pmatrix{ 1&-1& 0 \cr
                        0&1 &-1 \cr
                        0&0 &1  \cr}\ ,\ \ \bC {\bf C'}=\bI \ .
$$
Notice also that the matrix

$${\bf \tilde M}=\pmatrix{
-m_{1\da }e^{-{5\over\varepsilon}}&
0 &0 \cr
m_{21}e^{-{3\over\varepsilon}} &-m_{21}e^{-{3\over\varepsilon}}
  & 0 \cr 0 &
 m_{32}e^{-{1\over\varepsilon}} &-m_{32}e^{-{1\over\varepsilon}}
  \cr } \  $$
demonstrates the same in leading order eigenvalues asymptotic and the same
ultimate values of eigenvectors' components.

Although matrices $\bM $ and ${\bf \tilde M}$ outwardly are rather
different, however the solutions of the evolution equation

$$
{d\vec p\over dt}={\bf L}\vec p
$$
with guiding matrix ${\bf L}$ equal to $\bM ^*$ or to ${\bf \tilde M}^*$
are the same under exponentially large times $t=\tau e^{V/\varepsilon}$,
$V>0$, and in "slow" time" $\tau$ with exponential scale $V$ the
corresponding evolution $\vec p_V(\tau)=\lim\limits_{\varepsilon\to 0} 
\vec p(\tau e^{V/\varepsilon})$ has the form

$$\vec p_V(\tau)= {\bf C'}\left[ \lim\limits_{\varepsilon\to 0} {\rm
diag}\left\{ e^{\Lambda_1\tau\exp({V-5\over\varepsilon})},
e^{\Lambda_2\tau\exp({V-3\over\varepsilon})},
e^{\Lambda_1\tau\exp({V-1\over\varepsilon})}\right\}\right]\bC \vec
p_V(0)\ .
$$
Under exponential scale $V=1$ part of the initial distribution 
concentrated at the third state crosses to the second one. Such crossing 
rapidness at "slow" time $\tau$ is determined by the quantity 
$\Lambda_3=-m_{32}$.  Under exponential scale $V=3$ if $\tau\to\infty$ 
entire distribution turns out to be at the first state, corresponding 
rapidness is determined by the quantity $\Lambda_2=-m_{21}$. Finally at 
the scale $V=5$ killing of the process takes place and it is governed by 
the quantity $\Lambda_1=-m_{1\da }$. Note, that in spite of the killing 
probability $M_{3\da }=m_{3\da }e^{-4/\varepsilon}$ at state 3 is
greater than the killing probability $M_{1\da }=m_{1\da }e^{-5/
\varepsilon}$  at state 1, but  the part of distribution concentrated 
at state 3  crosses to the second state (and then to state 1) before 
killing at state 3 could take place.

Considered example shows from one side small matrix elements are not
negligible comparatively with the large ones. From another side
nevertheless one can neglect by some of the elements, may be not small (
here these are elements $M_{12}$, $M_{2\da}$, $M_{23}$, $M_{31}$,
$M_{3\da }$).  Which of elements do not affect on the spectrum and
on the evolution is determined by signless formulas (\ref{char}) and
(\ref{components}).

\begin{center}
\section{Discussion}
\end{center}

Obtained expressions for characteristic polynomial (\ref{char}) and
components of eigenvectors (\ref{components}) at the situation of arbitrary
matrices  give few substantial for practical computation in comparison
with usual methods. Moreover, the formula (\ref{circuit}) being only
reformulation of standard characteristic polynomial expression is
more preferable than (\ref{char}), because the volume of calculations
needed for determination of coefficients at powers of $\lambda$ using
(\ref{circuit}) is less than using (\ref{char}). The reason is that $G$
has less subgraphs being circuits than $G^\da $ has forests. In such
case signless formulas obtained are not of interest for practical
calculations.

However, the situation strongly changes in the case of large matrices
with non-negative elements which coefficients depend on small
(large) parameter and may possess on it different orders (such matrices,
as it was noted before, often appear in singular and stochastic problems).
Usual methods in such a situation are practically unsuitable due to large 
amount of calculations and low accuracy. On the contrary, formulas 
(\ref{char}) and (\ref{components}) due to sign absence at sums become 
extremely effective and allow to keep an eye on only higher order in 
asymptotic. This circumstance permits not only to carry out practical 
calculations, but also (which is highly essentially from theoretical point 
oh view) to determine singular limits and analyze structure spectrum 
properties. Here it is necessary to note that though for matrix $\bG $ 
with non-negative elements introducing quantities $g_{i\da 
}=-\sum\limits_{j=1}^Ng_{ij}$ are not positive however at the expense of 
spectrum shift it is always possible to make them non-negative ones.  The 
problem of picking out at characteristic polynomial terms having highly 
order is reducing to determination forests $F$ with maximal productivity 
$\pi_F$ and is a separate difficult one.  This problem is solved, the 
effective technique for determination of extreme forests is elaborated 
\cite{V5}. Analogical technique can be used for analysis of large very 
sparse matrices.

{\rm This work was  supported RFBR, grants N-99-01-00696 and
N-98-01-01063.} \vskip.5cm


\end{document}